\documentclass[review]{elsarticle}
\usepackage[utf8]{inputenc}
\usepackage{fbb}
\usepackage[margin = 1.5 cm]{geometry}
\usepackage{graphicx}
\usepackage{epstopdf}
\usepackage{float}
\usepackage{amsmath,amsfonts,amssymb,amsthm}
\usepackage{thmtools}
\usepackage{algorithm,algpseudocode}
\usepackage{xcolor}

\declaretheorem[style = plain,name = Proposition,parent = section,qed = \textnormal{\textemdash}]{proposicao}
\declaretheorem[style = plain,name = Definition,parent = section,qed = $\square $ ]{definicao}
\declaretheorem[style = plain,name = Theorem,parent = section,qed = \textnormal{\textemdash}]{teorema}

\journal{Journal of \LaTeX\ Templates}

\begin{document}

\begin{frontmatter}

\title{Piecewise Polynomial Interpolation Function Approach for Solving Nonlinear Programming Problems with Disjoint Feasible Regions: Mathematical Proofs}

\author[aluno]{Ricardo B. N. M. Pinheiro\corref{correspondingauthor}}
\cortext[correspondingauthor]{Corresponding Author}
\ead{ribenopi@hotmail.com}
\author[MAT]{Antonio R. Balbo\corref{colaborador}} \cortext[colaborador]{Colaborator}
\ead{antonio.balbo@unesp.br}
\author[UNESP]{Leonardo Nepomuceno\corref{supervisor}}
\cortext[supervisor]{Supervisor}
\ead{leonardo.nepomuceno@unesp.br}

\address[aluno]{Post-graduate program in Electrical Engineering, in the Engineering Faculty at Univ Estadual Paulista - UNESP, Bauru–SP, Brazil}
\address[UNESP]{Department of Electrical Engineering in the Engineering Faculty at Univ Estadual Paulista - UNESP, Bauru–SP, Brazil}
\address[MAT]{Department of Mathematics, Science Faculty at Univ Estadual Paulista - UNESP, Bauru–SP, Brazil}

\begin{abstract}
The Piecewise Polynomial Interpolation (PPI) function approach is aimed at solving  nonlinear programming  problems with disjoint feasible regions. In such problems, disjointedness  is generally  associated with prohibited operating zones, which correspond to bands of values that a variable is not allowed to assume. An analytical implication of such prohibited operating zones is to make the objective function, as well as  its domain, discontinuous. The PPI function  approach  consists in  replacing the constraints associated with prohibited operating zones by  an equivalent set of  equality and inequality constraints, thereby allowing the application of any efficient  gradient-based optimization method for solving the equivalent problem. In this paper, we present the definition of the PPI function and provide the mathematical proofs for its properties.
\end{abstract}


\end{frontmatter}


\section{Introduction}

A general Nonlinear Programming  Problem with Disjoint Feasible Regions (NLPDFR) \cite{smith_constraint_2013} is formulated as:
\begin{subequations}\label{eq1}
	\begin{align}
	\underset{\mathbf{x}\in {{\mathbb{R}}^{{{n}_{1}}}},\mathbf{p}\in {{\mathbb{R}}^{{{n}_{2}}}}}{\mathop{Min}}\, & f\left( \mathbf{x}\mathbf{,p} \right)  \\
	s.t: & \mathbf{g}\left( \mathbf{x}\mathbf{,p} \right)=\mathbf{0}  \\
	{} & \mathbf{h}\left( \mathbf{x}\mathbf{,p} \right)\le \mathbf{0}  \\
	{} & {{p}_{k}}\in \cup _{i=1}^{\text{NP}_{k}+1}\left[ p_{k,i}^{\min },p_{k,i}^{\max } \right],\forall k=1,...,{{n}_{2}},\label{eq1d} 
	\end{align}
\end{subequations}
where $f:{{\mathbb{R}}^{{{n}_{1}}}}\times {{\mathbb{R}}^{{{n}_{2}}}}\to \mathbb{R}$; $\mathbf{g}:{{\mathbb{R}}^{{{n}_{1}}}}\times {{\mathbb{R}}^{{{n}_{2}}}}\to {{\mathbb{R}}^{m}}$ and  $\mathbf{h}:{{\mathbb{R}}^{{{n}_{1}}}}\times {{\mathbb{R}}^{{{n}_{2}}}}\to {{\mathbb{R}}^{r}}$ are continuously differentiable linear or nonlinear functions, $\text{NP}_{k}$ is the number of Prohibited Operating Zones (POZ) of the variable ${{p}_{k}}$. We assume that $p_{k}^{\min }=p_{k,1}^{\min }<p_{k,1}^{\max }<...<p_{k,i}^{\min }<p_{k,i}^{\max }<...<p_{k,\text{NP}_{k}+1}^{\min }<p_{k,\text{NP}_{k}+1}^{\max }=p_{k}^{\max }$.

The constraints associated with the POZ are given in \eqref{eq1d}. Such constraints state that the variable  $p_{k}$ must stay inside the set composed by the union of the  $\text{NP}_{k} +1$ allowed operating zones. 

A practical example of a NLPDFR is the Economic Dispatch problem with Prohibited Operating Zones (EDPOZ) \cite{ding_bi-level_2015}. In the EDPOZ, the  objective is to calculate the generation dispatch for thermal units which minimize the fuel costs, while meeting the system demand and losses, enforcing the generation limits and avoiding the prohibited operating zones. In EDPOZ the POZ are associated with the operating zones where the thermal units are not allowed to operate due to mechanical operational problems, such as shaft bearing vibration.

\section{The Piecewise Polynomial Interpolation Approach}\label{secPOZ}

We start this section by considering the POZ constraints \eqref{eq1d}. Since $p_{k}^{\min }=p_{k,1}^{\min }<p_{k,1}^{\max }<...<p_{k,i}^{\min }<p_{k,i}^{\max }<...<p_{k,\text{NP}_{k}+1}^{\min }<p_{k,\text{NP}_{k}+1}^{\max }=p_{k}^{\max }$, it follows trivially that $\left[ p_{k,u}^{\min },p_{k,u}^{\max } \right]\cap \left[ p_{k,v}^{\min },p_{k,v}^{\max } \right] $ is the empty set, for any two distinct indexes $u$ and $v$. Moreover, for each index $i=1,..., \text{NP}_{k}+1$ and $k=1,...,n_{2}$, we have that all the elements in $\left[ p_{k,i}^{\min },p_{k,i}^{\max } \right]$ can be mapped by the one-to-one function (line segment) $d_{k,i}:\left[ 0,1 \right]\to \left[ p_{k,i}^{\min },p_{k,i}^{\max } \right]$, given by:
\begin{equation}\label{seg}
{d_{k,i}}\left( {{\alpha }_{k,i}} \right)={{\alpha }_{k,i}}p_{k,i}^{\max }+\left( 1-{{\alpha }_{k,i}} \right)p_{k,i}^{\min }.
\end{equation}

Based on \eqref{seg}, we establish  at the  Proposition \ref{prop1} an equivalent form of expressing POZ constraints \eqref{eq1d} in terms of ${d_{k,i}}$.
\begin{proposicao}\label{prop1}
	Let  $k\in {{\mathcal{G}}_{\text{POZ}}}$ and  ${{\mathcal{D}}_{k}}=\cup _{i=1}^{\text{NP}_{k}+1}\left\{ {{d}_{k,i}} \right\}$ be a discrete set of functions ${{d}_{k,i}}$ given in \eqref{seg}.   We say that ${{p}_{k}}\in \cup _{i=1}^{\text{NP}_{k}+1}\left[ p_{k,i}^{\min },p_{k,i}^{\max } \right]$ if, and only if, ${{d}_{k,i}}\left( {{\bar{\alpha } }_{k,i}} \right)={{p}_{k}}$ for some index $i$  and ${{\bar{\alpha }}_{k,i}}\in \left[ 0,1 \right]$.
\end{proposicao}

We present the Piecewise Polynomial Interpolation (PPI) function at  Definition \eqref{def2}. 
\begin{definicao}[PPI function]\label{def2}
	Let $k\in {{\mathcal{G}}_{\text{POZ}}}$ and $\mathbf{\alpha^{k}}=\left( {{\alpha }_{k,1}},...,{{\alpha }_{k,i}},...,{{\alpha }_{k,\text{NP}_{k}+1}} \right)$ be a vector in $\mathbb{I}=\times _{i=1}^{\text{NP}_{k}+1}\left[ 0,1 \right]\subset {{\mathbb{R}}^{\text{NP}_{k}+1}}$  and the function ${{d}_{k,i}}={{d}_{k,i}}\left( {{\alpha }_{k,i}} \right)$  presented in \eqref{seg}. The Piecewise Polynomial Interpolation (PPI)  is the function $m:\mathbb{R}\times \mathbb{I}\to \mathbb{R}$ given by:
	\begin{equation} \label{PPIsum}
	\begin{aligned}
	& m\left( {{p}_{k}},{{\mathbf{\alpha }}^{\mathbf{k}}} \right)=-\frac{\left( {{p}_{k}}-{{d}_{k,2}} \right)\left( {{p}_{k}}-{{d}_{k,1}} \right)}{{{d}_{k,2}}-{{d}_{k,1}}}+ \\ 
	& \ \ \ \ \ \ \ \ \ \ \ \ \ \ \ \ \ \ \ \ \ \ \ \ \ \ \ \ \ \ \ \sum\limits_{i=2}^{\text{NP}_{k}}{{{\left( -1 \right)}^{i}}\left[ \frac{\left( {{p}_{k}}-{{d}_{k,i}} \right)\left( {{p}_{k}}-{{d}_{k,i-1}} \right)}{{{d}_{k,i}}-{{d}_{k,i-1}}}+\frac{\left( {{p}_{k}}-{{d}_{k,i+1}} \right)\left( {{p}_{k}}-{{d}_{k,i}} \right)}{{{d}_{k,i+1}}-{{d}_{k,i}}} \right]u\left( {{p}_{k}},{{d}_{k,i}} \right)} 
	\end{aligned}
	\end{equation}
	where: 
	\begin{equation}\label{heaviside}
u\left( {{p}_{k}},{{d}_{k,i}} \right)=\left\{ \begin{array}{*{35}{l}}
1 & \text{if}\ {{p}_{k}}-{{d}_{k,i}}\ge 0  \\
0 & \text{otherwise}  \\
\end{array} \right.
\end{equation}
	is the Heaviside function.
\end{definicao}

The main properties of the PPI function are presented at   Theorem \ref{teo1}. The proofs for such theorem are provided in section \ref{prova}.

	\begin{teorema}\label{teo1}
	The PPI function $m\left( p_{k},\mathbf{\alpha^{k}}\right)$  \eqref{PPIsum} has the following properties:
	\begin{enumerate}
		\item  $m\left( p_{k},\mathbf{\alpha^{k}}\right)=0 \Leftrightarrow p_{k} \in \mathcal{D}_{k}$ and $m\left( {{p}_{k}},{{\mathbf{\alpha }}^{\mathbf{k}}} \right)\ne 0$, otherwise.\label{A1}		
		\item $\left| {{\nabla }_{{{p}_{k}}}}m\left( {{p}_{k}},{{\mathbf{\alpha }}^{\mathbf{k}}} \right) \right|=1\Leftrightarrow {{p}_{k}}\in {{\mathcal{D}}_{k}}$  (Normalization).\label{A2}		
		\item $m\left( p_{k},\mathbf{\alpha^{k}}\right)$ is of class $C^{1}$ (i.e., differentiable with  continuous partial derivatives).\label{A3} \qedhere
	\end{enumerate}
\end{teorema}

Using the properties described at the Theorem  \ref{teo1},  the proposed PPI approach consists in  replacing the POZ constraints  \eqref{eq1d} by the equivalent set of equality and inequality constraints given by:
\begin{subequations}\label{eqp7}
	\begin{align}
	{} & m\left( {{p}_{k}},{{\mathbf{\alpha }}^{\mathbf{k}}} \right)=0,\ \ \ \forall k=1,...,n_{2}\label{eqp7b}\\
	{} & 0\le {{\alpha }_{k,i}}\le 1,\,\,\,\,\,\,\,\,\  \ \,\forall k=1,...,n_{2}, i=1,...,\text{NP}_{k}+1.\label{eqp7c}
	\end{align}
\end{subequations}

Thus, by means of the PPI approach, the equivalent nonlinear programming problem to the NLPDFR \eqref{eq1}  is given by:
\begin{subequations}\label{eq3}
	\begin{align}
	\underset{\mathbf{x}\in {{\mathbb{R}}^{{{n}_{1}}}},\mathbf{p}\in {{\mathbb{R}}^{{{n}_{2}}}}}{\mathop{Min}}\, & f\left( \mathbf{x}\mathbf{,p} \right)  \\
	s.t: & \mathbf{g}\left( \mathbf{x}\mathbf{,p} \right)=\mathbf{0}  \\
	{} & \mathbf{h}\left( \mathbf{x}\mathbf{,p} \right)\le \mathbf{0}  \\
	{} & m\left( {{p}_{k}},{{\mathbf{\alpha }}^{\mathbf{k}}} \right)=0,\forall k=1,...,{{n}_{2}}  \\
	{} & 0\le {{\alpha }_{k,i}}\le 1,\ \forall k=1,...,{{n}_{2}}, i=1,...,\text{NP}_{k}+1.
	\end{align}  
\end{subequations}

The PPI approach has the following characteristics: (i) it does not require a reformulation by means of a MINLP problem for solving the original NLPDFR;  (ii) it eliminates the combinatorial nature and the disjoint feasible regions associated with POZ of the original NLPDFR problem and (iii) it allows for obtaining the solution of problem \eqref{eq1} by solving the equivalent problem \eqref{eq3}, by means of any efficient gradient-based method.

\section{Proofs for Theorem \ref{teo1}} \label{prova}

\begin{proof}
	\textbf{Statement \ref{A1}}.
	\newline
	$\left( \to  \right)$ 	Assume that $m\left( {{p}_{k}},{{\mathbf{\alpha }}^{\mathbf{k}}} \right)=0$.
	\newline
	
	Given a value for  ${{p}_{k}}$, we have that both the functions $u\left( {{p}_{k}},{{d}_{k,i}} \right)$ and  $m\left( {{p}_{k}},{{\mathbf{\alpha }}^{\mathbf{k}}} \right)$ assure that only one of the quadratic terms $\frac{\left( {{p}_{k}}-{{d}_{k,i}} \right)\left( {{p}_{k}}-{{d}_{k,i-1}} \right)}{{{d}_{k,i}}-{{d}_{k,i-1}}}$ or $\frac{\left( {{p}_{k}}-{{d}_{k,i+1}} \right)\left( {{p}_{k}}-{{d}_{k,i}} \right)}{{{d}_{k,i+1}}-{{d}_{k,i}}}$ will be null for some index  $i=1,...,\text{NP}_{k}$.	  Suppose that  this is verified for index $j$. Since ${{d}_{k,j}}-{{d}_{k,j-1}}\ne 0$ and ${{d}_{k,j+1}}-{{d}_{k,j}}\ne 0$ for all ${{\alpha }_{k,j-1}}$, ${{\alpha }_{k,j}}$, ${{\alpha }_{k,j+1}}$ $\in \left[ 0,1 \right]$, it is necessary that either ${{p}_{k}}-{{d}_{k,j-1}}=0$ or ${{p}_{k}}-{{d}_{k,j}}=0$ or ${{p}_{k}}-{{d}_{k,j+1}}=0$.	If ${{p}_{k}}-{{d}_{k,j}}=0$ (an analogous reasoning is used for the cases  ${{p}_{k}}-{{d}_{k,j-1}}=0$ and ${{p}_{k}}-{{d}_{k,j+1}}=0$), so, it follows from \eqref{seg} that there is only one  ${{\bar{\alpha }}_{k,j}}\in \left[ 0,1 \right]$ such that  ${{p}_{k}}-{{d}_{k,j}}\left( {{{\bar{\alpha }}}_{k,j}} \right)=0$. Therefore ${{p}_{k}}\in {{\mathcal{D}}_{k}}$.
	\\
	\newline
	$\left( \leftarrow  \right)$ Assume that ${{p}_{k}}\in {{\mathcal{D}}_{k}}$. 
	\newline
	
	So, for some index  $i=1,...,\text{NP}_{k}+1$, say $j$, there is only one  ${{\bar{\alpha }}_{k,j}}\in \left[ 0,1 \right]$  such that ${{p}_{k}}={{d}_{k,j}}\left( {{{\bar{\alpha }}}_{k,j}} \right)$ and ${{p}_{k}}-{{d}_{k,i}}\ne 0$ for all $i\ne j$ and ${{\alpha }_{k,i}}\in \left[ 0,1 \right]$. Hence, it follows from \eqref{heaviside} that the PPI function \eqref{PPIsum} becomes:
	\begin{align*}
	m\left( {{p}_{k}},{{\mathbf{\alpha }}^{\mathbf{k}}} \right)=\left\{ \begin{aligned}
	{{\left( -1 \right)}^{j}}\frac{\left( {{p}_{k}}-{{d}_{k,j+1}} \right)\left( {{p}_{k}}-{{d}_{k,j}}\left( {{{\bar{\alpha }}}_{k,j}} \right) \right)}{{{d}_{k,j+1}}-{{d}_{k,j}}\left( {{{\bar{\alpha }}}_{k,j}} \right)} &\ \  \text{, if}\ 1\le j\le \text{NP}_{k}  \\
	{{\left( -1 \right)}^{j-1}}\frac{\left( {{p}_{k}}-{{d}_{k,j}}\left( {{{\bar{\alpha }}}_{k,j}} \right) \right)\left( {{p}_{k}}-{{d}_{k,j-1}} \right)}{{{d}_{k,j}}\left( {{{\bar{\alpha }}}_{k,j}} \right)-{{d}_{k,j-1}}} &\ \  \text{, if}\ j=\text{NP}_{k}+1  \\
	\end{aligned} \right.
	\end{align*}
	
	Thus, it is immediate that $m\left( {{p}_{k}},{{\mathbf{\alpha }}^{\mathbf{k}}} \right)=0$.
	\\
	\newline
	\textbf{Statement 2}
	
	It follows from PPI function \eqref{PPIsum} that:
	\begin{equation*}
	{{\nabla }_{{{p}_{k}}}}m\left( {{p}_{k}},{{\mathbf{\alpha }}^{\mathbf{k}}} \right)=-\frac{2{{p}_{k}}-{{d}_{k,1}}-{{d}_{k,2}}}{{{d}_{k,2}}-{{d}_{k,1}}}+\sum\limits_{i=2}^{\text{NP}_{k}}{{{\left( -1 \right)}^{i}}\left[ \frac{2{{p}_{k}}-{{d}_{k,i}}-{{d}_{k,i-1}}}{{{d}_{k,i}}-{{d}_{k,i-1}}}+\frac{2{{p}_{k}}-{{d}_{k,i}}-{{d}_{k,i+1}}}{{{d}_{k,i+1}}-{{d}_{k,i}}} \right]u\left( {{p}_{k}},{{d}_{k,i}} \right)}.
	\end{equation*}
	
	It follows from Statement 1 that  $m\left( {{p}_{k}},{{\mathbf{\alpha }}^{\mathbf{k}}} \right)=0\Leftrightarrow {{p}_{k}}\in {{\mathcal{D}}_{k}}$. Thus, if  ${{p}_{k}}\in {{\mathcal{D}}_{k}}$ for some index $i=1,...,\text{NP}_{k}+1$, say $j$, there is only one ${{\bar{\alpha }}_{k,j}}\in \left[ 0,1 \right]$ such that  ${{p}_{k}}={{d}_{k,j}}\left( {{{\bar{\alpha }}}_{k,j}} \right)$ and ${{p}_{k}}-{{d}_{k,i}}\ne 0$ for all  $i\ne j$ and ${{\alpha }_{k,i}}\in \left[ 0,1 \right]$.  Hence, it follows from \eqref{heaviside} that the function  ${{\nabla }_{{{p}_{k}}}}m\left( {{p}_{k}},{{\mathbf{\alpha }}^{\mathbf{k}}} \right)$ becomes:
	\begin{equation*}
	{{\nabla }_{{{p}_{k}}}}m\left( {{p}_{k}},{{\mathbf{\alpha }}^{\mathbf{k}}} \right)=\left\{ \begin{aligned}
	{{\left( -1 \right)}^{j}}\frac{2{{p}_{k}}-{{d}_{k,j}}\left( {{{\bar{\alpha }}}_{k,j}} \right)-{{d}_{k,j+1}}}{{{d}_{k,j+1}}-{{d}_{k,j}}\left( {{{\bar{\alpha }}}_{k,j}} \right)} &\ \  \text{, if}\ 1\le j\le \text{NP}_{k}  \\
	{{\left( -1 \right)}^{j-1}}\frac{2{{p}_{k}}-{{d}_{k,j-1}}-{{d}_{k,j}}\left( {{{\bar{\alpha }}}_{k,j}} \right)}{{{d}_{k,j}}\left( {{{\bar{\alpha }}}_{k,j}} \right)-{{d}_{k,j-1}}} &\ \  \text{, if}\ j=\text{NP}_{k}+1  \\
	\end{aligned} \right.
	\end{equation*}
	
	Since ${{p}_{k}}={{d}_{k,j}}\left( {{{\bar{\alpha }}}_{k,j}} \right)$, it follows that  $\left| {{\nabla }_{{{p}_{k}}}}m\left( {{p}_{k}},{{\mathbf{\alpha }}^{\mathbf{k}}} \right) \right|=1$.
	\\
	\newline
	\textbf{Statement 3} 
	\\
	\newline
	The PPI function  $m\left( {{p}_{k}},{{\mathbf{\alpha }}^{\mathbf{k}}} \right)$ has partial derivatives at all points  $\left( {{p}_{k}},{{\mathbf{\alpha }}^{\mathbf{k}}} \right)\in \mathbb{R}\times _{i=1}^{\text{NP}_{k}+1}\left[ 0,1 \right]$. They are given by:
	\begin{equation}\label{gradp}
	{{\nabla }_{{{p}_{k}}}}m\left( {{p}_{k}},{{\mathbf{\alpha }}^{\mathbf{k}}} \right)=-\frac{2{{p}_{k}}-{{d}_{k,1}}-{{d}_{k,2}}}{{{d}_{k,2}}-{{d}_{k,1}}}+\sum\limits_{i=2}^{\text{NP}_{k}}{{{\left( -1 \right)}^{i}}\left[ \frac{2{{p}_{k}}-{{d}_{k,i}}-{{d}_{k,i-1}}}{{{d}_{k,i}}-{{d}_{k,i-1}}}+\frac{2{{p}_{k}}-{{d}_{k,i}}-{{d}_{k,i+1}}}{{{d}_{k,i+1}}-{{d}_{k,i}}} \right]u\left( {{p}_{k}},{{d}_{k,i}} \right)},
	\end{equation}
	
	\begin{equation}\label{grada1}
	{{\nabla }_{{{\alpha }_{k,1}}}}m\left( {{p}_{k}},{{\mathbf{\alpha }}^{\mathbf{k}}} \right)=-{{\left[ \frac{{{p}_{k}}-{{d}_{k,2}}}{{{d}_{k,2}}-{{d}_{k,1}}} \right]}^{2}}{{\Delta }_{k,1}}+{{\left[ \frac{{{p}_{k}}-{{d}_{k,2}}}{{{d}_{k,2}}-{{d}_{k,1}}} \right]}^{2}}{{\Delta }_{k,1}}u\left( {{p}_{k}},{{d}_{k,2}} \right),
	\end{equation}
	
	\begin{equation}\label{gradai}
	\begin{aligned}
	& {{\nabla }_{{{\alpha }_{k,i}}}}m\left( {{p}_{k}},{{\mathbf{\alpha }}^{\mathbf{k}}} \right)=-{{\left( -1 \right)}^{i-1}}{{\left[ \frac{{{p}_{k}}-{{d}_{k,i-1}}}{{{d}_{k,i}}-{{d}_{k,i-1}}} \right]}^{2}}{{\Delta }_{k,i}}u\left( {{p}_{k}},{{d}_{k,i-1}} \right)+ \\ 
	& \ \ \ \ \ \ \ \ \ \left[ -{{\left( -1 \right)}^{i}}{{\left[ \frac{{{p}_{k}}-{{d}_{k,i-1}}}{{{d}_{k,i}}-{{d}_{k,i-1}}} \right]}^{2}}+{{\left( -1 \right)}^{i}}{{\left[ \frac{{{p}_{k}}-{{d}_{k,i+1}}}{{{d}_{k,i+1}}-{{d}_{k,i}}} \right]}^{2}} \right]{{\Delta }_{k,i}}u\left( {{p}_{k}},{{d}_{k,i}} \right)+ \\ 
	& \ \ \ \ \ \ \ \ \ \ \ \ \ \ \ \ \ \ \ \ \ \ \ \ \ \ \ \ \ \ \ \ \ \ \ \ \ \ \ \ \ \ \ \ \ \ \ \ \ {{\left( -1 \right)}^{i+1}}{{\left[ \frac{{{p}_{k}}-{{d}_{k,i+1}}}{{{d}_{k,i+1}}-{{d}_{k,i}}} \right]}^{2}}{{\Delta }_{k,i}}u\left( {{p}_{k}},{{d}_{k,i+1}} \right)\text{, if  }1<i\le \text{NP}_{k} 
	\end{aligned}
	\end{equation}
	and
	\begin{equation}\label{gradaNP}
	{{\nabla }_{{{\alpha }_{k,\text{NP}_{k}+1}}}}m\left( {{p}_{k}},{{\mathbf{\alpha }}^{\mathbf{k}}} \right)=-{{\left( -1 \right)}^{\text{NP}_{k}}}{{\left[ \frac{{{p}_{k}}-{{d}_{k,\text{NP}_{k}}}}{{{d}_{k,\text{NP}_{k}+1}}-{{d}_{k,\text{NP}_{k}}}} \right]}^{2}}{{\Delta }_{k,\text{NP}_{k}+1}}u\left( {{p}_{k}},{{d}_{k,\text{NP}_{k}}} \right),
	\end{equation}
	where 
	\begin{equation}
		{{\Delta }_{k,i}}=p_{k,i}^{\max }-p_{k,i}^{\min },
	\end{equation}
	with $i=1,...,\text{NP}_{k}+1$.
	
	A sufficient condition for the PPI function \eqref{PPIsum} to be of class ${{C}^{1}}$ is to show that the partial derivatives ${{\nabla }_{{{p}_{k}}}}m\left( {{p}_{k}},{{\mathbf{\alpha }}^{\mathbf{k}}} \right)$ and ${{\nabla }_{{{\alpha }_{k,i}}}}m\left( {{p}_{k}},{{\mathbf{\alpha }}^{\mathbf{k}}} \right)$ are continuous in $\mathbb{R}\times _{i=1}^{\text{NP}_{k}+1}\left[ 0,1 \right]$. For this, we consider the continuity analysis of these partial derivatives at the following points:
	
	\begin{itemize} 
		\item $P=\left( {{{\bar{p}}}_{k}},{{{\mathbf{\bar{\alpha }}}}^{\mathbf{k}}} \right)$ such that  ${{\bar{p}}_{k}}-{{d}_{k,j}}\left( {{{\bar{\alpha }}}_{k,j}} \right)\ne 0$, for all  $j=1,...,\text{NP}_{k}+1$ and ${{\bar{\alpha }}_{k,j}}\in \left[ 0,1 \right]$ and;
		\item $Q=\left( {{{\bar{p}}}_{k}},{{{\mathbf{\bar{\alpha }}}}^{\mathbf{k}}} \right)$ such that  ${{\bar{p}}_{k}}-{{d}_{k,j}}\left( {{{\bar{\alpha }}}_{k,j}} \right)=0$ for some $j=1,...,\text{NP}_{k}+1$ and ${{\bar{\alpha }}_{k,j}}\in \left[ 0,1 \right]$.
	\end{itemize}
	
	When we consider the point  $P$, it follows from Heaviside funcion \eqref{heaviside},   ${{d}_{k,j+1}}-{{d}_{k,j}}\ne 0$ and/or ${{d}_{k,j}}-{{d}_{k,j-1}}\ne 0$ that both the functions  ${{\nabla }_{{{p}_{k}}}}m\left( {{p}_{k}},{{\mathbf{\alpha }}^{\mathbf{k}}} \right)$ and ${{\nabla }_{{{\alpha }_{k,i}}}}m\left( {{p}_{k}},{{\mathbf{\alpha }}^{\mathbf{k}}} \right)$ are rationals  (polynomial quotient) and, therefore, continuous at  $P$.
	
	In order to analyze the continuity of the functions ${{\nabla }_{{{p}_{k}}}}m\left( {{p}_{k}},{{\mathbf{\alpha }}^{\mathbf{k}}} \right)$ and  ${{\nabla }_{{{\alpha }_{k,i}}}}m\left( {{p}_{k}},{{\mathbf{\alpha }}^{\mathbf{k}}} \right)$ at the point $Q$,  we will consider the cases  $j=1$, $1<j\le \text{NP}_{k}$ and $j=\text{NP}_{k}+1$ separately.
	\\
	\newline
	\textbf{Continuity of the function ${{\nabla }_{{{p}_{k}}}}m\left( {{p}_{k}},{{\mathbf{\alpha }}^{\mathbf{k}}} \right)$}
	\\
	\newline
	Suppose $j=1$.
	\\
	\newline
	In an open ball  $\left\| \left( {{p}_{k}},{{\mathbf{\alpha }}^{\mathbf{k}}} \right)-Q \right\|<{{\delta }_{0}}$ with ${{\delta }_{0}}>0$ small enough, we have that the function  ${{\nabla }_{{{p}_{k}}}}m\left( {{p}_{k}},{{\mathbf{\alpha }}^{\mathbf{k}}} \right)$ is given by:
	\begin{equation*}
	{{\nabla }_{{{p}_{k}}}}m\left( {{p}_{k}},{{\mathbf{\alpha }}^{\mathbf{k}}} \right)=-\frac{2{{p}_{k}}-{{d}_{k,1}}-{{d}_{k,2}}}{{{d}_{k,2}}-{{d}_{k,1}}}.
	\end{equation*}
	It follows immediately that  ${{\nabla }_{{{p}_{k}}}}m\left( Q \right)=1$. Moreover, we have:
	\begin{equation*}
	\begin{aligned}
	\underset{\left( {{p}_{k}},{{\mathbf{\alpha }}^{\mathbf{k}}} \right)\to Q}{\mathop{\lim }}\,{{\nabla }_{{{p}_{k}}}}m\left( {{p}_{k}},{{\mathbf{\alpha }}^{\mathbf{k}}} \right)= & \underset{\left( {{p}_{k}},{{\mathbf{\alpha }}^{\mathbf{k}}} \right)\to Q}{\mathop{\lim }}\,-\frac{2{{p}_{k}}-{{d}_{k,1}}-{{d}_{k,2}}}{{{d}_{k,2}}-{{d}_{k,1}}}  \\
	{} & \underset{{{d}_{k,2}}-{{d}_{k,1}}\ne 0}{\mathop{=}}\,\frac{\underset{\left( {{p}_{k}},{{\mathbf{\alpha }}^{\mathbf{k}}} \right)\to Q}{\mathop{\lim }}\,\left( 2{{p}_{k}}-{{d}_{k,1}}-{{d}_{k,2}} \right)}{\underset{\left( {{p}_{k}},{{\mathbf{\alpha }}^{\mathbf{k}}} \right)\to Q}{\mathop{\lim }}\,\left( {{d}_{k,2}}-{{d}_{k,1}} \right)}  \\
	{} & =\frac{2{{{\bar{p}}}_{k}}-{{d}_{k,1}}\left( {{{\bar{\alpha }}}_{k,1}} \right)-d\left( {{\alpha }_{k,2}} \right)}{{{d}_{k,2}}\left( {{{\bar{\alpha }}}_{k,2}} \right)-d\left( {{{\bar{\alpha }}}_{k,1}} \right)}  \\
	{} & \underset{{{{\bar{p}}}_{k}}={{d}_{k,1}}\left( {{{\bar{\alpha }}}_{k,1}} \right)}{\mathop{=}}\,1.
	\end{aligned}
	\end{equation*}
	
	Since  $\underset{\left( {{p}_{k}},{{\mathbf{\alpha }}^{\mathbf{k}}} \right)\to Q}{\mathop{\lim }}\,{{\nabla }_{{{p}_{k}}}}m\left( {{p}_{k}},{{\mathbf{\alpha }}^{\mathbf{k}}} \right)=1={{\nabla }_{{{p}_{k}}}}m\left( Q \right)$, we conclude that   ${{\nabla }_{{{p}_{k}}}}m\left( {{p}_{k}},{{\mathbf{\alpha }}^{\mathbf{k}}} \right)$ is continuous at  $Q$ for $j=1$.
	\\
	\newline
	Suppose  $1<j\le \text{NP}_{k}$.  
	\\
	\newline
	In a open ball  $\left\| \left( {{p}_{k}},{{\mathbf{\alpha }}^{\mathbf{k}}} \right)-Q \right\|<{{\delta }_{1}}$ with ${{\delta }_{1}}>0$ small enough, we have that the function  ${{\nabla }_{{{p}_{k}}}}m\left( {{p}_{k}},{{\mathbf{\alpha }}^{\mathbf{k}}} \right)$ is given by:
	\begin{equation*}
	{{\nabla }_{{{p}_{k}}}}m\left( {{p}_{k}},{{\mathbf{\alpha }}^{\mathbf{k}}} \right)={{\left( -1 \right)}^{j-1}}\frac{2{{p}_{k}}-{{d}_{k,j}}-{{d}_{k,j-1}}}{{{d}_{k,j}}-{{d}_{k,j-1}}}+{{\left( -1 \right)}^{j}}\left[ \frac{2{{p}_{k}}-{{d}_{k,j}}-{{d}_{k,j-1}}}{{{d}_{k,j}}-{{d}_{k,j-1}}}+\frac{2{{p}_{k}}-{{d}_{k,j}}-{{d}_{k,j+1}}}{{{d}_{k,j+1}}-{{d}_{k,j}}} \right]u\left( {{p}_{k}},{{d}_{k,j}} \right).
	\end{equation*}
	
	It follows immediately that  ${{\nabla }_{{{p}_{k}}}}m\left( Q \right)=-{{\left( -1 \right)}^{j}}$. Let us verify the behavior of function   ${{\nabla }_{{{p}_{k}}}}m\left( {{p}_{k}},{{\mathbf{\alpha }}^{\mathbf{k}}} \right)$ around  $Q$ when ${{p}_{k}}$ tends to the right and left of ${{d}_{k,j}}\left( {{{\bar{\alpha }}}_{k,j}} \right)$.
	
	Suppose  ${{p}_{k}}<{{d}_{k,j}}$. Thus, it follows that  $u\left( {{p}_{k}},{{d}_{k,j}} \right)=0$ and
	\begin{equation*}
	\begin{aligned}
	\underset{\left( {{p}_{k}},{{\mathbf{\alpha }}^{\mathbf{k}}} \right)\to Q}{\mathop{\lim }}\,{{\nabla }_{{{p}_{k}}}}m\left( {{p}_{k}},{{\mathbf{\alpha }}^{\mathbf{k}}} \right) & =\underset{\left( {{p}_{k}},{{\mathbf{\alpha }}^{\mathbf{k}}} \right)\to Q}{\mathop{\lim }}\,{{\left( -1 \right)}^{j-1}}\frac{2{{p}_{k}}-{{d}_{k,j}}-{{d}_{k,j-1}}}{{{d}_{k,j}}-{{d}_{k,j-1}}}  \\
	{} & \underset{{{d}_{k,j}}\ne {{d}_{k,j-1}}}{\mathop{=}}\,{{\left( -1 \right)}^{j-1}}\frac{\underset{\left( {{p}_{k}},{{\mathbf{\alpha }}^{\mathbf{k}}} \right)\to Q}{\mathop{\lim }}\,\left[ 2{{p}_{k}}-{{d}_{k,j}}-{{d}_{k,j-1}} \right]}{\underset{\left( {{p}_{k}},{{\mathbf{\alpha }}^{\mathbf{k}}} \right)\to Q}{\mathop{\lim }}\,\left[ {{d}_{k,j}}-{{d}_{k,j-1}} \right]}  \\
	{} & ={{\left( -1 \right)}^{j-1}}\frac{2{{{\bar{p}}}_{k}}-{{d}_{k,j}}\left( {{{\bar{\alpha }}}_{k,j}} \right)-{{d}_{k,j-1}}\left( {{{\bar{\alpha }}}_{k,j-1}} \right)}{{{d}_{k,j}}\left( {{{\bar{\alpha }}}_{k,j}} \right)-{{d}_{k,j-1}}\left( {{{\bar{\alpha }}}_{k,j-1}} \right)}  \\
	{} & \underset{{{{\bar{p}}}_{k}}={{d}_{k,j}}\left( {{{\bar{\alpha }}}_{k,j}} \right)}{\mathop{=}}\,-{{\left( -1 \right)}^{j}}. 
	\end{aligned}
	\end{equation*}
	
	Suppose ${{p}_{k}}>{{d}_{k,j}}$. Thus, it follows that  $u\left( {{p}_{k}},{{d}_{k,j}} \right)=1$ and
	\begin{equation*}
	\begin{aligned}
	\underset{\left( {{p}_{k}},{{\mathbf{\alpha }}^{\mathbf{k}}} \right)\to Q}{\mathop{\lim }}\,{{\nabla }_{{{p}_{k}}}}m\left( {{p}_{k}},{{\mathbf{\alpha }}^{\mathbf{k}}} \right) & =\underset{\left( {{p}_{k}},{{\mathbf{\alpha }}^{\mathbf{k}}} \right)\to Q}{\mathop{\lim }}\,{{\left( -1 \right)}^{j}}\frac{2{{p}_{k}}-{{d}_{k,j}}-{{d}_{k,j+1}}}{{{d}_{k,j+1}}-{{d}_{k,j}}}  \\
	{} & \underset{{{d}_{k,j+1}}\ne {{d}_{k,j}}}{\mathop{=}}\,{{\left( -1 \right)}^{j}}\frac{\underset{\left( {{p}_{k}},{{\mathbf{\alpha }}^{\mathbf{k}}} \right)\to Q}{\mathop{\lim }}\,\left[ 2{{p}_{k}}-{{d}_{k,j}}-{{d}_{k,j+1}} \right]}{\underset{\left( {{p}_{k}},{{\mathbf{\alpha }}^{\mathbf{k}}} \right)\to Q}{\mathop{\lim }}\,\left[ {{d}_{k,j+1}}-{{d}_{k,j}} \right]}  \\
	{} & ={{\left( -1 \right)}^{j}}\frac{2{{{\bar{p}}}_{k}}-{{d}_{k,j}}\left( {{{\bar{\alpha }}}_{k,j}} \right)-{{d}_{k,j+1}}\left( {{{\bar{\alpha }}}_{k,j+1}} \right)}{{{d}_{k,j+1}}\left( {{{\bar{\alpha }}}_{k,j+1}} \right)-{{d}_{k,j}}\left( {{{\bar{\alpha }}}_{k,j}} \right)}  \\
	{} & \underset{{{{\bar{p}}}_{k}}={{d}_{k,j}}\left( {{{\bar{\alpha }}}_{k,j}} \right)}{\mathop{=}}\,-{{\left( -1 \right)}^{j}}.  
	\end{aligned}
	\end{equation*}
	
	Since $\underset{\left( {{p}_{k}},{{\mathbf{\alpha }}^{\mathbf{k}}} \right)\to Q}{\mathop{\lim }}\,{{\nabla }_{{{p}_{k}}}}m\left( {{p}_{k}},{{\mathbf{\alpha }}^{\mathbf{k}}} \right)=-{{\left( -1 \right)}^{j}}={{\nabla }_{{{p}_{k}}}}m\left( Q \right)$, we conclude that ${{\nabla }_{{{p}_{k}}}}m\left( {{p}_{k}},{{\mathbf{\alpha }}^{\mathbf{k}}} \right)$ is continuous at  $Q$ for  $1<j\le \text{NP}_{k}$.
	\\
	\newline
	Suppose  $j=\text{NP}_{k}+1$.  
	\\
	\newline
	In an open ball  $\left\| \left( {{p}_{k}},{{\mathbf{\alpha }}^{\mathbf{k}}} \right)-Q \right\|<{{\delta }_{2}}$ with ${{\delta }_{2}}>0$ small enough, we have that the function ${{\nabla }_{{{p}_{k}}}}m\left( {{p}_{k}},{{\mathbf{\alpha }}^{\mathbf{k}}} \right)$ is given by:
	\begin{equation*}
	{{\nabla }_{{{p}_{k}}}}m\left( {{p}_{k}},{{\mathbf{\alpha }}^{\mathbf{k}}} \right)={{\left( -1 \right)}^{\text{NP}_{k}}}\frac{2{{p}_{k}}-{{d}_{k,\text{NP}_{k}}}-{{d}_{k,\text{NP}_{k}+1}}}{{{d}_{k,\text{NP}_{k}+1}}-{{d}_{k,\text{NP}_{k}}}}.
	\end{equation*}
	It follows immediately that ${{\nabla }_{{{p}_{k}}}}m\left( {{p}_{k}},{{\mathbf{\alpha }}^{\mathbf{k}}} \right)={{\left( -1 \right)}^{\text{NP}_{k}}}$. Moreover, we have:
	\begin{equation*}
	\begin{aligned}
	\underset{\left( {{p}_{k}},{{\mathbf{\alpha }}^{\mathbf{k}}} \right)\to Q}{\mathop{\lim }}\,{{\nabla }_{{{p}_{k}}}}m\left( {{p}_{k}},{{\mathbf{\alpha }}^{\mathbf{k}}} \right) & =\underset{\left( {{p}_{k}},{{\mathbf{\alpha }}^{\mathbf{k}}} \right)\to Q}{\mathop{\lim }}\,{{\left( -1 \right)}^{\text{NP}_{k}}}\frac{2{{p}_{k}}-{{d}_{k,\text{NP}_{k}}}-{{d}_{k,\text{NP}_{k}+1}}}{{{d}_{k,\text{NP}_{k}+1}}-{{d}_{k,\text{NP}_{k}}}}  \\
	{} & ={{\left( -1 \right)}^{\text{NP}_{k}}}\frac{2{{{\bar{p}}}_{k}}-{{d}_{k,\text{NP}_{k}}}\left( {{{\bar{\alpha }}}_{k,\text{NP}_{k}}} \right)-{{d}_{k,\text{NP}_{k}+1}}\left( {{{\bar{\alpha }}}_{k,\text{NP}_{k}+1}} \right)}{{{d}_{k,\text{NP}_{k}+1}}\left( {{{\bar{\alpha }}}_{k,\text{NP}_{k}+1}} \right)-{{d}_{k,\text{NP}_{k}}}\left( {{{\bar{\alpha }}}_{k,\text{NP}_{k}}} \right)}  \\
	{} & \underset{{{{\bar{p}}}_{k}}={{d}_{k,\text{NP}_{k}+1}}\left( {{{\bar{\alpha }}}_{k,\text{NP}_{k}+1}} \right)}{\mathop{=}}\,{{\left( -1 \right)}^{\text{NP}_{k}}}.  
	\end{aligned}
	\end{equation*}
	Since $\underset{\left( {{p}_{k}},{{\mathbf{\alpha }}^{\mathbf{k}}} \right)\to Q}{\mathop{\lim }}\,{{\nabla }_{{{p}_{k}}}}m\left( {{p}_{k}},{{\mathbf{\alpha }}^{\mathbf{k}}} \right)={{\left( -1 \right)}^{j-1}}={{\nabla }_{{{p}_{k}}}}m\left( Q \right)$, we conclude that  ${{\nabla }_{{{p}_{k}}}}m\left( {{p}_{k}},{{\mathbf{\alpha }}^{\mathbf{k}}} \right)$ is continuous at $Q$ for  $j=\text{NP}_{k}+1$.
	
	Therefore, the function  ${{\nabla }_{{{p}_{k}}}}m\left( {{p}_{k}},{{\mathbf{\alpha }}^{\mathbf{k}}} \right)$ is continuous in $\mathbb{R}\times _{i=1}^{\text{NP}_{k}+1}\left[ 0,1 \right]$.
	\\
	\newline
	\textbf{Continuity of the function ${{\nabla }_{{{\alpha }_{k,i}}}}m\left( {{p}_{k}},{{\mathbf{\alpha }}^{\mathbf{k}}} \right)$}
	\\
	\newline
	Suppose $j=1$.
	\\
	\newline
	In an open ball  $\left\| \left( {{p}_{k}},{{\mathbf{\alpha }}^{\mathbf{k}}} \right)-Q \right\|<{{\delta }_{3}}$ with ${{\delta }_{3}}>0$ small enough, we have that the function  ${{\nabla }_{{{\alpha }_{k,i}}}}m\left( {{p}_{k}},{{\mathbf{\alpha }}^{\mathbf{k}}} \right)$ is given by:
	\begin{equation*}
	{{\nabla }_{{{\alpha }_{k,i}}}}m\left( {{p}_{k}},{{\mathbf{\alpha }}^{\mathbf{k}}} \right)=\left\{ \begin{array}{*{35}{l}}
	\begin{aligned}
	{{\left( -1 \right)}^{i}}{{\left[ \frac{{{p}_{k}}-{{d}_{k,i+1}}}{{{d}_{k,i+1}}-{{d}_{k,i}}} \right]}^{2}}{{\Delta }_{k,i}} 
	\end{aligned}& ,\text{if  }i=1  \\
	
	\begin{aligned}-{{\left( -1 \right)}^{i-1}}{{\left[ \frac{{{p}_{k}}-{{d}_{k,i-1}}}{{{d}_{k,i}}-{{d}_{k,i-1}}} \right]}^{2}}{{\Delta }_{k,i}}u\left( {{p}_{k}},{{d}_{k,i-1}} \right)
	\end{aligned} & ,\text{if  }i=2  \\
	0 & \text{, if  }2<i\le \text{NP}_{k}+1  \\
	\end{array} \right.
	\end{equation*}
	
	It follows immediately that ${{\nabla }_{{{\alpha }_{k,1}}}}m\left( Q \right)=-{{\Delta }_{k,1}}$ and ${{\nabla }_{{{\alpha }_{k,2}}}}m\left( Q \right)=0$. For $i=1$, we have:
	\begin{equation*}
	\begin{aligned}
	\underset{\left( {{p}_{k}},{{\mathbf{\alpha }}^{\mathbf{k}}} \right)\to Q}{\mathop{\lim }}\,{{\nabla }_{{{\alpha }_{k,1}}}}m\left( {{p}_{k}},{{\mathbf{\alpha }}^{\mathbf{k}}} \right) & =\underset{\left( {{p}_{k}},{{\mathbf{\alpha }}^{\mathbf{k}}} \right)\to Q}{\mathop{\lim }}\,-{{\left[ \frac{{{p}_{k}}-{{d}_{k,2}}}{{{d}_{k,2}}-{{d}_{k,1}}} \right]}^{2}}{{\Delta }_{k,1}}  \\
	{} & \underset{{{d}_{k,2}}\ne {{d}_{k,1}}}{\mathop{=}}\,-{{\left[ \frac{\underset{\left( {{p}_{k}},{{\mathbf{\alpha }}^{\mathbf{k}}} \right)\to Q}{\mathop{\lim }}\,\left( {{p}_{k}}-{{d}_{k,2}} \right)}{\underset{\left( {{p}_{k}},{{\mathbf{\alpha }}^{\mathbf{k}}} \right)\to Q}{\mathop{\lim }}\,\left( {{d}_{k,2}}-{{d}_{k,1}} \right)} \right]}^{2}}{{\Delta }_{k,1}}  \\
	{} & =-{{\left[ \frac{{{{\bar{p}}}_{k}}-{{d}_{k,2}}\left( {{{\bar{\alpha }}}_{k,2}} \right)}{{{d}_{k,2}}\left( {{{\bar{\alpha }}}_{k,2}} \right)-{{d}_{k,1}}\left( {{{\bar{\alpha }}}_{k,1}} \right)} \right]}^{2}}{{\Delta }_{k,1}}  \\
	{} & \underset{{{{\bar{p}}}_{k}}={{d}_{k,1}}\left( {{{\bar{\alpha }}}_{k,1}} \right)}{\mathop{=}}\,-{{\Delta }_{k,1}}. 
	\end{aligned}_{{}}^{{}}
	\end{equation*}
	
	Since $\underset{\left( {{p}_{k}},{{\mathbf{\alpha }}^{\mathbf{k}}} \right)\to Q}{\mathop{\lim }}\,{{\nabla }_{{{\alpha }_{k,1}}}}m\left( {{p}_{k}},{{\mathbf{\alpha }}^{\mathbf{k}}} \right)={{\nabla }_{{{\alpha }_{k,1}}}}m\left( Q \right)=-{{\Delta }_{k,1}}$, we conclude that  ${{\nabla }_{{{\alpha }_{k,1}}}}m\left( {{p}_{k}},{{\mathbf{\alpha }}^{\mathbf{k}}} \right)$ is continuous at $Q$ for $j=1$.
	
	For $i=2$, let us  verify the behavior of function  ${{\nabla }_{{{\alpha }_{k,2}}}}m\left( {{p}_{k}},{{\mathbf{\alpha }}^{\mathbf{k}}} \right)$ around $Q$ when  ${{p}_{k}}$ tends to the right and left of  ${{d}_{k,1}}\left( {{{\bar{\alpha }}}_{k,1}} \right)$.
	Suppose that ${{p}_{k}}<{{d}_{k,1}}$. Thus, it follows that $u\left( {{p}_{k}},{{d}_{k,j}} \right)=0$ and $\underset{\left( {{p}_{k}},{{\mathbf{\alpha }}^{\mathbf{k}}} \right)\to Q}{\mathop{\lim }}\,{{\nabla }_{{{\alpha }_{k,2}}}}m\left( {{p}_{k}},{{\mathbf{\alpha }}^{\mathbf{k}}} \right)=0$.
	Suppose that  ${{p}_{k}}>{{d}_{k,1}}$. Thus, it follows that $u\left( {{p}_{k}},{{d}_{k,j}} \right)=1$  and
	\begin{equation*}
	\begin{aligned}
	\underset{\left( {{p}_{k}},{{\mathbf{\alpha }}^{\mathbf{k}}} \right)\to Q}{\mathop{\lim }}\,{{\nabla }_{{{\alpha }_{k,2}}}}m\left( {{p}_{k}},{{\mathbf{\alpha }}^{\mathbf{k}}} \right) & =\underset{\left( {{p}_{k}},{{\mathbf{\alpha }}^{\mathbf{k}}} \right)\to Q}{\mathop{\lim }}\,{{\left[ \frac{{{p}_{k}}-{{d}_{k,1}}}{{{d}_{k,2}}-{{d}_{k,1}}} \right]}^{2}}{{\Delta }_{k,2}}  \\
	{} & \underset{{{d}_{k,2}}\ne {{d}_{k,1}}}{\mathop{=}}\,\frac{\underset{\left( {{p}_{k}},{{\mathbf{\alpha }}^{\mathbf{k}}} \right)\to Q}{\mathop{\lim }}\,{{\left[ {{p}_{k}}-{{d}_{k,1}} \right]}^{2}}}{\underset{\left( {{p}_{k}},{{\mathbf{\alpha }}^{\mathbf{k}}} \right)\to Q}{\mathop{\lim }}\,{{\left[ {{d}_{k,2}}-{{d}_{k,1}} \right]}^{2}}}{{\Delta }_{k,2}}  \\
	{} & ={{\left[ \frac{{{{\bar{p}}}_{k}}-{{d}_{k,1}}\left( {{{\bar{\alpha }}}_{k,1}} \right)}{{{d}_{k,2}}\left( {{{\bar{\alpha }}}_{k,2}} \right)-{{d}_{k,1}}\left( {{{\bar{\alpha }}}_{k,1}} \right)} \right]}^{2}}{{\Delta }_{k,2}}  \\
	{} & \underset{{{{\bar{p}}}_{k}}={{d}_{k,1}}\left( {{{\bar{\alpha }}}_{k,1}} \right)}{\mathop{=}}\,0.  
	\end{aligned}
	\end{equation*}
	
	Since $\underset{\left( {{p}_{k}},{{\mathbf{\alpha }}^{\mathbf{k}}} \right)\to Q}{\mathop{\lim }}\,{{\nabla }_{{{\alpha }_{k,2}}}}m\left( {{p}_{k}},{{\mathbf{\alpha }}^{\mathbf{k}}} \right)={{\nabla }_{{{\alpha }_{k,2}}}}m\left( Q \right)=0$, we conclude that  ${{\nabla }_{{{\alpha }_{k,2}}}}m\left( {{p}_{k}},{{\mathbf{\alpha }}^{\mathbf{k}}} \right)$ is continuous at  $Q$ for $j=1$.
	\\
	\newline
	Suppose $1<j\le \text{NP}_{k}$.
	\\
	\newline
	In an open ball $\left\| \left( {{p}_{k}},{{\mathbf{\alpha }}^{\mathbf{k}}} \right)-Q \right\|<{{\delta }_{3}}$ with ${{\delta }_{3}}>0$ small enough, we have that the function ${{\nabla }_{{{\alpha }_{k,i}}}}m\left( {{p}_{k}},{{\mathbf{\alpha }}^{\mathbf{k}}} \right)$ is given by:
	\begin{equation*}
	{{\nabla }_{{{\alpha }_{k,1}}}}m\left( {{p}_{k}},{{\mathbf{\alpha }}^{\mathbf{k}}} \right)=\left\{ \begin{array}{*{35}{l}}
	\begin{aligned}-{{\left[ \frac{{{p}_{k}}-{{d}_{k,2}}}{{{d}_{k,2}}-{{d}_{k,1}}} \right]}^{2}}{{\Delta }_{k,1}}+{{\left[ \frac{{{p}_{k}}-{{d}_{k,2}}}{{{d}_{k,2}}-{{d}_{k,1}}} \right]}^{2}}{{\Delta }_{k,1}}u\left( {{p}_{k}},{{d}_{k,2}} \right)
	\end{aligned} & ,\text{if}\ j=2  \\
	0 & ,\text{if}\ j>2  \\
	\end{array} \right.
	\end{equation*}
	
	\begin{equation*}
	{{\nabla }_{{{\alpha }_{k,i}}}}m\left( {{p}_{k}},{{\mathbf{\alpha }}^{\mathbf{k}}} \right)=\left\{ \begin{array}{*{35}{l}}
	0 & ,\ \text{if}\ j<i-1  \\
	\begin{aligned}-{{\left( -1 \right)}^{i-1}}{{\left[ \frac{{{p}_{k}}-{{d}_{k,i-1}}}{{{d}_{k,i}}-{{d}_{k,i-1}}} \right]}^{2}}{{\Delta }_{k,i}}u\left( {{p}_{k}},{{d}_{k,i-1}} \right)
	\end{aligned} & ,\text{if  }j=i-1  \\
	\begin{aligned}
	& -{{\left( -1 \right)}^{i-1}}{{\left[ \frac{{{p}_{k}}-{{d}_{k,i-1}}}{{{d}_{k,i}}-{{d}_{k,i-1}}} \right]}^{2}}{{\Delta }_{k,i}}- \\ 
	& \left[ {{\left( -1 \right)}^{i}}{{\left[ \frac{{{p}_{k}}-{{d}_{k,i-1}}}{{{d}_{k,i}}-{{d}_{k,i-1}}} \right]}^{2}}-{{\left( -1 \right)}^{i}}{{\left[ \frac{{{p}_{k}}-{{d}_{k,i+1}}}{{{d}_{k,i+1}}-{{d}_{k,i}}} \right]}^{2}} \right]{{\Delta }_{k,i}}u\left( {{p}_{k}},{{d}_{k,i}} \right) \\ 
	\end{aligned} & ,\text{if  }j=i  \\
	\begin{aligned}{{\left( -1 \right)}^{i}}{{\left[ \frac{{{p}_{k}}-{{d}_{k,i+1}}}{{{d}_{k,i+1}}-{{d}_{k,i}}} \right]}^{2}}{{\Delta }_{k,i}}+{{\left( -1 \right)}^{i+1}}{{\left[ \frac{{{p}_{k}}-{{d}_{k,i+1}}}{{{d}_{k,i+1}}-{{d}_{k,i}}} \right]}^{2}}{{\Delta }_{k,i}}u\left( {{p}_{k}},{{d}_{k,i+1}} \right)
	\end{aligned} & ,\text{if  }j=i+1  \\
	0 & \text{, if  }j>i+1  \\
	\end{array} \right.
	\end{equation*}
	$\text{if  }1<i\le \text{NP}_{k}$ and
	\begin{equation*}
	{{\nabla }_{{{\alpha }_{k,\text{NP}_{k}+1}}}}m\left( {{p}_{k}},{{\mathbf{\alpha }}^{\mathbf{k}}} \right)=\left\{ \begin{array}{*{35}{l}}
	0 & \text{, if  }j<\text{NP}_{k}  \\
	\begin{aligned}-{{\left( -1 \right)}^{\text{NP}_{k}}}{{\left[ \frac{{{p}_{k}}-{{d}_{k,\text{NP}_{k}}}}{{{d}_{k,\text{NP}_{k}+1}}-{{d}_{k,\text{NP}_{k}}}} \right]}^{2}}{{\Delta }_{k,\text{NP}_{k}+1}}u\left( {{p}_{k}},{{d}_{k,\text{NP}_{k}}} \right)
	\end{aligned} & \text{, if  }j=\text{NP}_{k}.  \\
	\end{array} \right.
	\end{equation*}
	
	It follows immediately that
	\begin{equation*}
	\begin{aligned}
	& {{\nabla }_{{{\alpha }_{k,1}}}}m\left( Q \right)=0 \\ 
	& {{\nabla }_{{{\alpha }_{k,i}}}}m\left( Q \right)=\left\{ \begin{array}{*{35}{l}}
	0 & ,\ \text{if}\ j<i-1 & {}  \\
	0 & ,\text{if  }j=i-1 & {}  \\
	{{\left( -1 \right)}^{i}}{{\Delta }_{k,i}} & ,\text{if  }j=i & ,\text{if  }1<i\le \text{NP}_{k}  \\
	0 & ,\text{if  }j=i+1 & {}  \\
	0  & \text{, if  }j>i+1 & {}  \\
	\end{array} \right. \\ 
	& {{\nabla }_{{{\alpha }_{k,\text{NP}_{k}+1}}}}m\left( Q \right)=0.
	\end{aligned}
	\end{equation*}
	
	For  $i=1$, let us to verify the behavior of the function   ${{\nabla }_{{{\alpha }_{k,1}}}}m\left( {{p}_{k}},{{\mathbf{\alpha }}^{\mathbf{k}}} \right)$ around $Q$ when ${{p}_{k}}$ tends to the right and left of  ${{d}_{k,2}}\left( {{{\bar{\alpha }}}_{k,2}} \right)$.
	Suppose that ${{p}_{k}}<{{d}_{k,2}}$. Thus, it follows that  $u\left( {{p}_{k}},{{d}_{k,j}} \right)=0$ and
	\begin{equation*}
	\begin{aligned}
	\underset{\left( {{p}_{k}},{{\mathbf{\alpha }}^{\mathbf{k}}} \right)\to Q}{\mathop{\lim }}\,{{\nabla }_{{{\alpha }_{k,1}}}}m\left( {{p}_{k}},{{\mathbf{\alpha }}^{\mathbf{k}}} \right) & =\underset{\left( {{p}_{k}},{{\mathbf{\alpha }}^{\mathbf{k}}} \right)\to Q}{\mathop{\lim }}\,-{{\left[ \frac{{{p}_{k}}-{{d}_{k,2}}}{{{d}_{k,2}}-{{d}_{k,1}}} \right]}^{2}}{{\Delta }_{k,1}}  \\
	{} & \underset{{{d}_{k,2}}\ne {{d}_{k,1}}}{\mathop{=}}\,-\frac{\underset{\left( {{p}_{k}},{{\mathbf{\alpha }}^{\mathbf{k}}} \right)\to Q}{\mathop{\lim }}\,{{\left[ {{p}_{k}}-{{d}_{k,2}} \right]}^{2}}}{\underset{\left( {{p}_{k}},{{\mathbf{\alpha }}^{\mathbf{k}}} \right)\to Q}{\mathop{\lim }}\,{{\left[ {{d}_{k,2}}-{{d}_{k,1}} \right]}^{2}}}{{\Delta }_{k,1}}  \\
	{} & =-{{\left[ \frac{{{{\bar{p}}}_{k}}-{{d}_{k,2}}\left( {{{\bar{\alpha }}}_{k,2}} \right)}{{{d}_{k,2}}\left( {{{\bar{\alpha }}}_{k,2}} \right)-{{d}_{k,1}}\left( {{{\bar{\alpha }}}_{k,1}} \right)} \right]}^{2}}{{\Delta }_{k,1}}  \\
	{} & \underset{{{{\bar{p}}}_{k}}={{d}_{k,2}}\left( {{{\bar{\alpha }}}_{k,2}} \right)}{\mathop{=}}\,0.  \\
	\end{aligned}.
	\end{equation*}
	
	Suppose that ${{p}_{k}}>{{d}_{k,2}}$. Thus, it follows that  $u\left( {{p}_{k}},{{d}_{k,j}} \right)=1$ and $\underset{\left( {{p}_{k}},{{\mathbf{\alpha }}^{\mathbf{k}}} \right)\to Q}{\mathop{\lim }}\,{{\nabla }_{{{\alpha }_{k,1}}}}m\left( {{p}_{k}},{{\mathbf{\alpha }}^{\mathbf{k}}} \right)=0$. 
	
	Since  $\underset{\left( {{p}_{k}},{{\mathbf{\alpha }}^{\mathbf{k}}} \right)\to Q}{\mathop{\lim }}\,{{\nabla }_{{{\alpha }_{k,1}}}}m\left( {{p}_{k}},{{\mathbf{\alpha }}^{\mathbf{k}}} \right)={{\nabla }_{{{\alpha }_{k,1}}}}m\left( Q \right)=0$, we conclude that ${{\nabla }_{{{\alpha }_{k,1}}}}m\left( {{p}_{k}},{{\mathbf{\alpha }}^{\mathbf{k}}} \right)$ is continuous at $Q$ for $1<j\le \text{NP}_{k}$.
	
	For $1<i\le \text{NP}_{k}$, let us to verify the behavior of the function  ${{\nabla }_{{{\alpha }_{k,i}}}}m\left( {{p}_{k}},{{\mathbf{\alpha }}^{\mathbf{k}}} \right)$ around $Q$ when ${{p}_{k}}$ tends to the right and left of ${{d}_{k,j-1}}\left( {{{\bar{\alpha }}}_{k,j-1}} \right)$,  ${{d}_{k,j}}\left( {{{\bar{\alpha }}}_{k,j}} \right)$ and ${{d}_{k,j+1}}\left( {{{\bar{\alpha }}}_{k,j+1}} \right)$.

	Suppose that ${{p}_{k}}<{{d}_{k,j-1}}$. Thus, it follows that  $u\left( {{p}_{k}},{{d}_{k,j-1}} \right)=0$ and $\underset{\left( {{p}_{k}},{{\mathbf{\alpha }}^{\mathbf{k}}} \right)\to Q}{\mathop{\lim }}\,{{\nabla }_{{{\alpha }_{k,i}}}}m\left( {{p}_{k}},{{\mathbf{\alpha }}^{\mathbf{k}}} \right)=0$.
	
	Suppose ${{p}_{k}}>{{d}_{k,j-1}}$. Thus, it follows that $u\left( {{p}_{k}},{{d}_{k,j-1}} \right)=1$ and
	\begin{equation*}
	\begin{aligned}
	\underset{\left( {{p}_{k}},{{\mathbf{\alpha }}^{\mathbf{k}}} \right)\to Q}{\mathop{\lim }}\,{{\nabla }_{{{\alpha }_{k,i}}}}m\left( {{p}_{k}},{{\mathbf{\alpha }}^{\mathbf{k}}} \right) & =\underset{\left( {{p}_{k}},{{\mathbf{\alpha }}^{\mathbf{k}}} \right)\to Q}{\mathop{\lim }}\,-{{\left( -1 \right)}^{i-1}}{{\left[ \frac{{{p}_{k}}-{{d}_{k,i-1}}}{{{d}_{k,i}}-{{d}_{k,i-1}}} \right]}^{2}}{{\Delta }_{k,i}}  \\
	{} & \underset{{{d}_{k,i}}\ne {{d}_{k,i-1}}}{\mathop{=}}\,-{{\left( -1 \right)}^{i-1}}\frac{\underset{\left( {{p}_{k}},{{\mathbf{\alpha }}^{\mathbf{k}}} \right)\to Q}{\mathop{\lim }}\,{{\left[ {{p}_{k}}-{{d}_{k,i-1}} \right]}^{2}}}{\underset{\left( {{p}_{k}},{{\mathbf{\alpha }}^{\mathbf{k}}} \right)\to Q}{\mathop{\lim }}\,{{\left[ {{d}_{k,i}}-{{d}_{k,i-1}} \right]}^{2}}}{{\Delta }_{k,i}}  \\
	{} & =-{{\left( -1 \right)}^{i-1}}{{\left[ \frac{{{{\bar{p}}}_{k}}-{{d}_{k,i-1}}\left( {{{\bar{\alpha }}}_{k,i-1}} \right)}{{{d}_{k,i}}\left( {{{\bar{\alpha }}}_{k,i}} \right)-{{d}_{k,i-1}}\left( {{{\bar{\alpha }}}_{k,i-1}} \right)} \right]}^{2}}{{\Delta }_{k,i}}  \\
	{} & \underset{{{{\bar{p}}}_{k}}={{d}_{k,i-1}}\left( {{{\bar{\alpha }}}_{k,i-1}} \right)}{\mathop{=}}\,0.  
	\end{aligned}
	\end{equation*}
	
	Since $\underset{\left( {{p}_{k}},{{\mathbf{\alpha }}^{\mathbf{k}}} \right)\to Q}{\mathop{\lim }}\,{{\nabla }_{{{\alpha }_{k,i}}}}m\left( {{p}_{k}},{{\mathbf{\alpha }}^{\mathbf{k}}} \right)={{\nabla }_{{{\alpha }_{k,i}}}}m\left( Q \right)=0$, we conclude that ${{\nabla }_{{{\alpha }_{k,i}}}}m\left( {{p}_{k}},{{\mathbf{\alpha }}^{\mathbf{k}}} \right)$ is continuous at  $Q$ for  $j=i-1$.  
	
	Suppose that ${{p}_{k}}<{{d}_{k,j}}$. Thus, it follows that  $u\left( {{p}_{k}},{{d}_{k,j}} \right)=0$ and 
	\begin{equation*}
	\begin{aligned}
	\underset{\left( {{p}_{k}},{{\mathbf{\alpha }}^{\mathbf{k}}} \right)\to Q}{\mathop{\lim }}\,{{\nabla }_{{{\alpha }_{k,i}}}}m\left( {{p}_{k}},{{\mathbf{\alpha }}^{\mathbf{k}}} \right) & =\underset{\left( {{p}_{k}},{{\mathbf{\alpha }}^{\mathbf{k}}} \right)\to Q}{\mathop{\lim }}\,-{{\left( -1 \right)}^{i-1}}{{\left[ \frac{{{p}_{k}}-{{d}_{k,i-1}}}{{{d}_{k,i}}-{{d}_{k,i-1}}} \right]}^{2}}{{\Delta }_{k,i}}  \\
	{} & \underset{{{d}_{k,i}}\ne {{d}_{k,i-1}}}{\mathop{=}}\,-{{\left( -1 \right)}^{i-1}}\frac{\underset{\left( {{p}_{k}},{{\mathbf{\alpha }}^{\mathbf{k}}} \right)\to Q}{\mathop{\lim }}\,{{\left[ {{p}_{k}}-{{d}_{k,i-1}} \right]}^{2}}}{\underset{\left( {{p}_{k}},{{\mathbf{\alpha }}^{\mathbf{k}}} \right)\to Q}{\mathop{\lim }}\,{{\left[ {{d}_{k,i}}-{{d}_{k,i-1}} \right]}^{2}}}{{\Delta }_{k,i}}  \\
	{} & =-{{\left( -1 \right)}^{i-1}}{{\left[ \frac{{{{\bar{p}}}_{k}}-{{d}_{k,i-1}}\left( {{{\bar{\alpha }}}_{k,i-1}} \right)}{{{d}_{k,i}}\left( {{{\bar{\alpha }}}_{k,i}} \right)-{{d}_{k,i-1}}\left( {{{\bar{\alpha }}}_{k,i-1}} \right)} \right]}^{2}}{{\Delta }_{k,i}}  \\
	{} & \underset{{{{\bar{p}}}_{k}}={{d}_{k,i}}\left( {{{\bar{\alpha }}}_{k,i}} \right)}{\mathop{=}}\,{{\left( -1 \right)}^{i}}{{\Delta }_{k,i}}.  
	\end{aligned}
	\end{equation*}
	
	Suppose that  ${{p}_{k}}>{{d}_{k,j}}$. Thus, it follows that $u\left( {{p}_{k}},{{d}_{k,j}} \right)=1$ and
	\begin{equation*}
	\begin{aligned}
	\underset{\left( {{p}_{k}},{{\mathbf{\alpha }}^{\mathbf{k}}} \right)\to Q}{\mathop{\lim }}\,{{\nabla }_{{{\alpha }_{k,i}}}}m\left( {{p}_{k}},{{\mathbf{\alpha }}^{\mathbf{k}}} \right) & =\underset{\left( {{p}_{k}},{{\mathbf{\alpha }}^{\mathbf{k}}} \right)\to Q}{\mathop{\lim }}\,{{\left( -1 \right)}^{i}}{{\left[ \frac{{{p}_{k}}-{{d}_{k,i+1}}}{{{d}_{k,i+1}}-{{d}_{k,i}}} \right]}^{2}}{{\Delta }_{k,i}}  \\
	{} & \underset{{{d}_{k,i+1}}\ne {{d}_{k,i}}}{\mathop{=}}\,{{\left( -1 \right)}^{i}}\frac{\underset{\left( {{p}_{k}},{{\mathbf{\alpha }}^{\mathbf{k}}} \right)\to Q}{\mathop{\lim }}\,{{\left[ {{p}_{k}}-{{d}_{k,i+1}} \right]}^{2}}}{\underset{\left( {{p}_{k}},{{\mathbf{\alpha }}^{\mathbf{k}}} \right)\to Q}{\mathop{\lim }}\,{{\left[ {{d}_{k,i+1}}-{{d}_{k,i}} \right]}^{2}}}{{\Delta }_{k,i}}  \\
	{} & ={{\left( -1 \right)}^{i}}{{\left[ \frac{{{{\bar{p}}}_{k}}-{{d}_{k,i+1}}\left( {{{\bar{\alpha }}}_{k,i+1}} \right)}{{{d}_{k,i+1}}\left( {{{\bar{\alpha }}}_{k,i+1}} \right)-{{d}_{k,i}}\left( {{{\bar{\alpha }}}_{k,i}} \right)} \right]}^{2}}{{\Delta }_{k,i}}  \\
	{} & \underset{{{{\bar{p}}}_{k}}={{d}_{k,i}}\left( {{{\bar{\alpha }}}_{k,i}} \right)}{\mathop{=}}\,{{\left( -1 \right)}^{i}}{{\Delta }_{k,i}}.  
	\end{aligned}
	\end{equation*}
	
	Since  $\underset{\left( {{p}_{k}},{{\mathbf{\alpha }}^{\mathbf{k}}} \right)\to Q}{\mathop{\lim }}\,{{\nabla }_{{{\alpha }_{k,i}}}}m\left( {{p}_{k}},{{\mathbf{\alpha }}^{\mathbf{k}}} \right)={{\nabla }_{{{\alpha }_{k,i}}}}m\left( Q \right)={{\left( -1 \right)}^{i}}{{\Delta }_{k,i}}$, we conclude that ${{\nabla }_{{{\alpha }_{k,i}}}}m\left( {{p}_{k}},{{\mathbf{\alpha }}^{\mathbf{k}}} \right)$ is continuous at  $Q$ for $j=i$.  
	
	Suppose that  ${{p}_{k}}<{{d}_{k,j+1}}$. Thus, it follows that  $u\left( {{p}_{k}},{{d}_{k,j+1}} \right)=0$ and
	\begin{equation*}
	\begin{aligned}
	\underset{\left( {{p}_{k}},{{\mathbf{\alpha }}^{\mathbf{k}}} \right)\to Q}{\mathop{\lim }}\,{{\nabla }_{{{\alpha }_{k,i}}}}m\left( {{p}_{k}},{{\mathbf{\alpha }}^{\mathbf{k}}} \right) & =\underset{\left( {{p}_{k}},{{\mathbf{\alpha }}^{\mathbf{k}}} \right)\to Q}{\mathop{\lim }}\,{{\left( -1 \right)}^{i}}{{\left[ \frac{{{p}_{k}}-{{d}_{k,i+1}}}{{{d}_{k,i+1}}-{{d}_{k,i}}} \right]}^{2}}{{\Delta }_{k,i}}  \\
	{} & \underset{{{d}_{k,i+1}}\ne {{d}_{k,i}}}{\mathop{=}}\,{{\left( -1 \right)}^{i}}\frac{\underset{\left( {{p}_{k}},{{\mathbf{\alpha }}^{\mathbf{k}}} \right)\to Q}{\mathop{\lim }}\,{{\left[ {{p}_{k}}-{{d}_{k,i+1}} \right]}^{2}}}{\underset{\left( {{p}_{k}},{{\mathbf{\alpha }}^{\mathbf{k}}} \right)\to Q}{\mathop{\lim }}\,{{\left[ {{d}_{k,i+1}}-{{d}_{k,i}} \right]}^{2}}}{{\Delta }_{k,i}}  \\
	{} & ={{\left( -1 \right)}^{i}}{{\left[ \frac{{{{\bar{p}}}_{k}}-{{d}_{k,i+1}}\left( {{{\bar{\alpha }}}_{k,i+1}} \right)}{{{d}_{k,i+1}}\left( {{{\bar{\alpha }}}_{k,i+1}} \right)-{{d}_{k,i}}\left( {{{\bar{\alpha }}}_{k,i}} \right)} \right]}^{2}}{{\Delta }_{k,i}}  \\
	{} & \underset{{{{\bar{p}}}_{k}}={{d}_{k,i+1}}\left( {{{\bar{\alpha }}}_{k,i+1}} \right)}{\mathop{=}}\,0.  
	\end{aligned}
	\end{equation*} 
	
	Suppose that ${{p}_{k}}>{{d}_{k,j+1}}$. Thus, it follows that $u\left( {{p}_{k}},{{d}_{k,j+1}} \right)=1$ and $\underset{\left( {{p}_{k}},{{\mathbf{\alpha }}^{\mathbf{k}}} \right)\to Q}{\mathop{\lim }}\,{{\nabla }_{{{\alpha }_{k,i}}}}m\left( {{p}_{k}},{{\mathbf{\alpha }}^{\mathbf{k}}} \right)=0$. 
	
	Since  $\underset{\left( {{p}_{k}},{{\mathbf{\alpha }}^{\mathbf{k}}} \right)\to Q}{\mathop{\lim }}\,{{\nabla }_{{{\alpha }_{k,i}}}}m\left( {{p}_{k}},{{\mathbf{\alpha }}^{\mathbf{k}}} \right)={{\nabla }_{{{\alpha }_{k,i}}}}m\left( Q \right)=0$, we conclude that  ${{\nabla }_{{{\alpha }_{k,i}}}}m\left( {{p}_{k}},{{\mathbf{\alpha }}^{\mathbf{k}}} \right)$ is continuous at $Q$ for $j=i+1$.  
	
	For $i=\text{NP}_{k}+1$, let us verify the behavior of the function  ${{\nabla }_{{{\alpha }_{k,\text{NP}_{k}+1}}}}m\left( {{p}_{k}},{{\mathbf{\alpha }}^{\mathbf{k}}} \right)$ around $Q$ when ${{p}_{k}}$ tends to the right and left of ${{d}_{k,\text{NP}_{k}}}\left( {{{\bar{\alpha }}}_{k,\text{NP}_{k}}} \right)$.
	
	Suppose that ${{p}_{k}}<{{d}_{k,\text{NP}_{k}}}$. Thus, it follows that  $u\left( {{p}_{k}},{{d}_{k,\text{NP}_{k}}} \right)=0$ and  $\underset{\left( {{p}_{k}},{{\mathbf{\alpha }}^{\mathbf{k}}} \right)\to Q}{\mathop{\lim }}\,{{\nabla }_{{{\alpha }_{k,\text{NP}_{k}+1}}}}m\left( {{p}_{k}},{{\mathbf{\alpha }}^{\mathbf{k}}} \right)=0$.
	
	Suppose that  ${{p}_{k}}>{{d}_{k,\text{NP}_{k}}}$. Thus, it follows that  $u\left( {{p}_{k}},{{d}_{k,\text{NP}_{k}}} \right)=1$ and
	\begin{equation*}
	\begin{aligned}
	\underset{\left( {{p}_{k}},{{\mathbf{\alpha }}^{\mathbf{k}}} \right)\to Q}{\mathop{\lim }}\,{{\nabla }_{{{\alpha }_{k,\text{NP}_{k}+1}}}}m\left( {{p}_{k}},{{\mathbf{\alpha }}^{\mathbf{k}}} \right) & =\underset{\left( {{p}_{k}},{{\mathbf{\alpha }}^{\mathbf{k}}} \right)\to Q}{\mathop{\lim }}\,-{{\left( -1 \right)}^{\text{NP}_{k}}}{{\left[ \frac{{{p}_{k}}-{{d}_{k,\text{NP}_{k}}}}{{{d}_{k,\text{NP}_{k}+1}}-{{d}_{k,\text{NP}_{k}}}} \right]}^{2}}{{\Delta }_{k,\text{NP}_{k}+1}}  \\
	{} & \underset{{{d}_{k,i}}\ne {{d}_{k,i-1}}}{\mathop{=}}\,-{{\left( -1 \right)}^{\text{NP}_{k}}}\frac{\underset{\left( {{p}_{k}},{{\mathbf{\alpha }}^{\mathbf{k}}} \right)\to Q}{\mathop{\lim }}\,{{\left[ {{p}_{k}}-{{d}_{k,\text{NP}_{k}}} \right]}^{2}}}{\underset{\left( {{p}_{k}},{{\mathbf{\alpha }}^{\mathbf{k}}} \right)\to Q}{\mathop{\lim }}\,{{\left[ {{d}_{k,\text{NP}_{k}+1}}-{{d}_{k,\text{NP}_{k}}} \right]}^{2}}}{{\Delta }_{k,\text{NP}_{k}+1}}  \\
	{} & =-{{\left( -1 \right)}^{\text{NP}_{k}}}{{\left[ \frac{{{{\bar{p}}}_{k}}-{{d}_{k,\text{NP}_{k}}}\left( {{{\bar{\alpha }}}_{k,\text{NP}_{k}}} \right)}{{{d}_{k,\text{NP}_{k}+1}}\left( {{{\bar{\alpha }}}_{k,\text{NP}_{k}+1}} \right)-{{d}_{k,\text{NP}_{k}}}\left( {{{\bar{\alpha }}}_{k,\text{NP}_{k}}} \right)} \right]}^{2}}{{\Delta }_{k,i}}  \\
	{} & \underset{{{{\bar{p}}}_{k}}={{d}_{k,\text{NP}_{k}}}\left( {{{\bar{\alpha }}}_{k,\text{NP}_{k}}} \right)}{\mathop{=}}\,0.  
	\end{aligned}
	\end{equation*}
	
	Since $\underset{\left( {{p}_{k}},{{\mathbf{\alpha }}^{\mathbf{k}}} \right)\to Q}{\mathop{\lim }}\,{{\nabla }_{{{\alpha }_{k,\text{NP}_{k}+1}}}}m\left( {{p}_{k}},{{\mathbf{\alpha }}^{\mathbf{k}}} \right)={{\nabla }_{{{\alpha }_{k,\text{NP}_{k}+1}}}}m\left( Q \right)=0$, we conclude that ${{\nabla }_{{{\alpha }_{k,i}}}}m\left( {{p}_{k}},{{\mathbf{\alpha }}^{\mathbf{k}}} \right)$ is continuous at $Q$ for $j=\text{NP}_{k}$.
	\\
	\newline
	Finally, suppose $j= \text{NP}_{k}+1$.
	\\
	\newline
	In an open ball  $\left\| \left( {{p}_{k}},{{\mathbf{\alpha }}^{\mathbf{k}}} \right)-Q \right\|<{{\delta }_{4}}$ with ${{\delta }_{4}}>0$ small enough, we have that the function  ${{\nabla }_{{{\alpha }_{k,i}}}}m\left( {{p}_{k}},{{\mathbf{\alpha }}^{\mathbf{k}}} \right)$ is given by:
	\begin{equation*}
	{{\nabla }_{{{\alpha }_{k,i}}}}m\left( {{p}_{k}},{{\mathbf{\alpha }}^{\mathbf{k}}} \right)=\left\{ \begin{array}{*{35}{l}}
	0 & \text{, if  }i=1,...,\text{NP}_{k}  \\
	\begin{aligned}-{{\left( -1 \right)}^{i-1}}{{\left[ \frac{{{p}_{k}}-{{d}_{k,i-1}}}{{{d}_{k,i}}-{{d}_{k,i-1}}} \right]}^{2}}{{\Delta }_{k,i}} 
	\end{aligned}& ,\text{if  }i=\text{NP}_{k}+1 
	\end{array} \right.
	\end{equation*}
	
	It follows immediately that 
	\begin{equation*}
	{{\nabla }_{{{\alpha }_{k,i}}}}m\left( Q \right)=\left\{ \begin{array}{*{35}{l}}
	0 & \text{, if  }1\le i\le \text{NP}_{k}  \\
	\begin{aligned}-{{\left( -1 \right)}^{i-1}}{{\Delta }_{k,i}}
	\end{aligned} & ,\text{if  }i=\text{NP}_{k}+1  
	\end{array} \right.
	\end{equation*}
	
	For $1\le i\le \text{NP}_{k}$, the continuity of ${{\nabla }_{{{\alpha }_{k,i}}}}m\left( {{p}_{k}},{{\mathbf{\alpha }}^{\mathbf{k}}} \right)$ at $Q$ when $j=\text{NP}_{k}+1$ is verified immediately. For $i=\text{NP}_{k}+1$, we have:
	\begin{equation*}
	\begin{aligned}
	\underset{\left( {{p}_{k}},{{\mathbf{\alpha }}^{\mathbf{k}}} \right)\to Q}{\mathop{\lim }}\,{{\nabla }_{{{\alpha }_{k,\text{NP}_{k}+1}}}}m\left( {{p}_{k}},{{\mathbf{\alpha }}^{\mathbf{k}}} \right) & =\underset{\left( {{p}_{k}},{{\mathbf{\alpha }}^{\mathbf{k}}} \right)\to Q}{\mathop{\lim }}\,-{{\left( -1 \right)}^{\text{NP}_{k}}}{{\left[ \frac{{{p}_{k}}-{{d}_{k,\text{NP}_{k}}}}{{{d}_{k,\text{NP}_{k}+1}}-{{d}_{k,\text{NP}_{k}}}} \right]}^{2}}{{\Delta }_{k,\text{NP}_{k}+1}}  \\
	{} & \underset{{{d}_{k,\text{NP}_{k}+}}\ne {{d}_{k,\text{NP}_{k}}}}{\mathop{=}}\,-{{\left( -1 \right)}^{\text{NP}_{k}}}{{\left[ \frac{\underset{\left( {{p}_{k}},{{\mathbf{\alpha }}^{\mathbf{k}}} \right)\to Q}{\mathop{\lim }}\,\left( {{p}_{k}}-{{d}_{k,\text{NP}_{k}}} \right)}{\underset{\left( {{p}_{k}},{{\mathbf{\alpha }}^{\mathbf{k}}} \right)\to Q}{\mathop{\lim }}\,\left( {{d}_{k,\text{NP}_{k}+1}}-{{d}_{k,\text{NP}_{k}}} \right)} \right]}^{2}}{{\Delta }_{k,\text{NP}_{k}+1}}  \\
	{} & =-{{\left( -1 \right)}^{\text{NP}_{k}}}{{\left[ \frac{{{{\bar{p}}}_{k}}-{{d}_{k,\text{NP}_{k}}}\left( {{{\bar{\alpha }}}_{k,\text{NP}_{k}}} \right)}{{{d}_{k,\text{NP}_{k}+1}}\left( {{{\bar{\alpha }}}_{k,\text{NP}_{k}+1}} \right)-{{d}_{k,\text{NP}_{k}}}\left( {{{\bar{\alpha }}}_{k,\text{NP}_{k}}} \right)} \right]}^{2}}{{\Delta }_{k,\text{NP}_{k}+1}}  \\
	{} & \underset{{{{\bar{p}}}_{k}}={{d}_{k,\text{NP}_{k}+1}}\left( {{{\bar{\alpha }}}_{k,\text{NP}_{k}+1}} \right)}{\mathop{=}}\,-{{\left( -1 \right)}^{\text{NP}_{k}}}{{\Delta }_{k,\text{NP}_{k}+1}}.  
	\end{aligned}
	\end{equation*}
	
	Since $\underset{\left( {{p}_{k}},{{\mathbf{\alpha }}^{\mathbf{k}}} \right)\to Q}{\mathop{\lim }}\,{{\nabla }_{{{\alpha }_{k,\text{NP}_{k}+1}}}}m\left( {{p}_{k}},{{\mathbf{\alpha }}^{\mathbf{k}}} \right)={{\nabla }_{{{\alpha }_{k,\text{NP}_{k}+1}}}}m\left( Q \right)=-{{\left( -1 \right)}^{\text{NP}_{k}}}{{\Delta }_{k,\text{NP}_{k}+1}}$, we conclude that  ${{\nabla }_{{{\alpha }_{k,\text{NP}_{k}+1}}}}m\left( {{p}_{k}},{{\mathbf{\alpha }}^{\mathbf{k}}} \right)$ is continuous $Q$ for $j=\text{NP}_{k}+1$.

	Therefore, the function  ${{\nabla }_{{{\alpha }_{i}}}}m\left( {{p}_{k}},{{\mathbf{\alpha }}^{\mathbf{k}}} \right)$ é is continuous in $\mathbb{R}\times _{i=1}^{\text{NP}_{k}+1}\left[ 0,1 \right]$.
	
	We conclude that the function $m\left( {{p}_{k}},{{\mathbf{\alpha }}^{\mathbf{k}}} \right)$ is of class ${{C}^{1}}$ (differentiable with partial derivatives  ${{\nabla }_{{{p}_{k}}}}m\left( {{p}_{k}},{{\mathbf{\alpha }}^{\mathbf{k}}} \right)$ and ${{\nabla }_{{{\alpha }_{i}}}}m\left( {{p}_{k}},{{\mathbf{\alpha }}^{\mathbf{k}}} \right)$ continuous). The proof is complete.
\end{proof}

%
%
%

\bibliography{references}

\end{document}